# Comment

**Grace Wahba**

## 1. INTRODUCTION

The authors are to be commended for jumping in to describe support vector machines (SVMs), not an easy thing to do since the the literature for SVMs has grown at least exponentially in the last few years. A Google search for "support vector machines" gave "about 1,180,000" hits as of this writing. The authors have nevertheless made a nice selection of important points to emphasize. As noted, SVMs were proposed for classification in the early 1990s by arguments like those behind Figure 1 in their paper. The use of SVMs grew rapidly among computer scientists, as it was found that they worked very well in all kinds of practical applications. The theoretical underpinnings that went with the original proposals were different than those in the classical statistical literature, for example, those related to Bayes risk, and so had less impact in the statistical literature. The convergence of SVMs and regularization methods (or, rather the convergence of the "SVM community" and the "regularization community") was a major impetus in the study of the (classical) statistical properties of the SVM. One point at which this convergence took place was at an American Mathematical Society meeting at Mt. Holyoke in 1996. The speaker was describing the SVM with the so-called kernel trick when an anonymous person at the back of the room remarked that the SVM with the kernel trick was the solution to an optimization problem in a reproducing kernel Hilbert space (RKHS). Once it was clear to statisticians that the SVM can be obtained as the result of an optimization/regularization problem in a RKHS, tools known to statisticians in this context were rapidly employed to show how the SVM could be modified to take into account nonrepresentative sample sizes, unequal misclassification costs and more than two classes, and to show in each case that it directly targets the Bayes risk under very general circumstances (see also [5, 8]). Thus, a "classical" explanation of why they work so well was provided.

*Grace Wahba is the IJ Schoenberg-Hilldale Professor of Statistics, Department of Statistics, University of Wisconsin, 1300 University Avenue, Madison, Wisconsin 53706, USA e-mail: wahba@stat.wisc.edu, and is also a member of the Computer Sciences Department and the Biostatistics and Medical Informatics Department.*

## 2. MERCER'S KERNELS AND POSITIVE DEFINITE FUNCTIONS

Let $\mathcal{T}$ be a.d.o. (any dirty old) domain and let $K(s,t), s,t \in \mathcal{T}$, be a symmetric, positive definite function of two variables; $K$ is said to be positive definite if for any $n$, and any $t_1, \ldots, t_n \in \mathcal{T}$, the $n \times n$ matrix with $ij$th element $K(t_i, t_j)$ is nonnegative definite. In the early SVM literature, as well as in the present paper, the kernel is described as having a representation $K(s,t) = \sum_{\nu=1}^{\infty} \lambda_\nu \Phi_\nu(s) \Phi_\nu(t)$. Here the (nonnegative) $\lambda_\nu$ and the $\Phi_\nu$ are the eigenvalues and eigenvectors of $K$. A representation as in this sum is sufficient for $K$ to be positive definite (see [13] on the Mercer Hilbert–Schmidt theorem), but the so-called radial basis functions (RBF) popular in machine learning, of the form $K(s,t) = k(\|s-t\|)$, $s, t$ in Euclidean $d$-space $E^d$, do not have a countable sequence of eigenvalues and eigenvectors—complex exponentials play the role of eigenvectors (see [3]). The Gaussian kernel $K_c(x,y) = e^{-\|x-y\|^2/c}$ is such an example. Although the notion of a countable expansion was used in uncoupling the linear SVM from its linearity restriction (and seems to be repeated over and over), the lack of a countable set of eigenvectors and eigenvalues does not affect the use of the Gaussian kernel or any other positive definite function in an SVM; as the authors note, only values of $K$ are needed. The RBF probably just do not want to be called "Mercer's kernels" (!). Positive definite functions are sometimes called reproducing kernels, relating to their association with RKHS [1].

Given a collection of objects (which could be vectors, images, sounds, graphs, texts, trees, ...) in

a.d.o. domain $\mathcal{T}$, a positive definite matrix with $ij$ entry $K(i,j)$ defines a (squared) distance $d_{ij}$ between the $i$th and $j$th object as

$$d_{ij} = K(i,i) + K(j,j) - 2K(i,j)$$

(and, in addition, this distance comes with an inner product). It can be argued that using distance between objects, defined in some way, is truly fundamental to classification and, therefore, positive definite kernels, since they provide a distance, play a fundamental role.

## 3. LARGE MARGIN CLASSIFIERS AND REGULARIZATION

Referring to equation (3.1) in the main paper, note that the elementary cost function $(1 - y_i f(\mathbf{x}_i))_+$ depends only on $\tau_i = y_i f(\mathbf{x}_i)$. If $y_i$ and $f(\mathbf{x}_i)$ have the same sign, then $f$ will classify $y_i$ correctly, and if they have different signs, then $f$ will classify $y_i$ incorrectly. The term $\tau_i$ is frequently called the margin, and classifiers that depend on the data only through $\tau$ are called large margin classifiers. The cost function $c(\tau) = (\tau)_+$ is called the misclassification counter, and it would be considered the ideal cost function if it were not for the fact that it leads to a nonconvex, nontractable optimization problem. Considering Bernoulli data coded as $y_i = 1$ or $y_i = 0$, the penalized likelihood estimate, where the cost function is the negative log likelihood, goes back at least to [12]. In that paper, members of the exponential family were considered as cost functions and it was natural to put the log likelihood in the canonical form for distributions in the exponential family. Thus the log likelihood for Bernoulli data is parameterized by the logit $f(x) = \log p(x)/(1 - p(x))$. However, if Bernoulli data are recoded as $y_i = \pm 1$, then the log likelihood (cost function) becomes $\mathcal{L}(y,f) = (1 + e^{-yf})$. Since thresholding $p(x)$ at $p = 1/2$ is equivalent to thresholding at $f = 0$, the penalized log likelihood estimate is also a large margin classifier.

It turns out that there are lots of large margin classifiers with the property that the sign of the estimate that minimizes

$$\frac{1}{n} \sum_{i=1}^{n} c(y_i f(\mathbf{x}_i)) + \mu \|f\|_K^2$$

tends to the sign of the log odds ratio, assuming that the problem is tuned adequately and that the RKHS associated with $K$ is rich enough for the problem at hand. The following rather amazing result is from [6]: Let $c(z) < c(-z)$, every $z > 0$, and let $c'(0) \neq 0$ exist. If $Ec(Yf((X)|\mathbf{X} = \mathbf{x}))$ has a global minimizer $\bar{f}(\mathbf{x})$ and $f(\mathbf{x}) \neq 0$, then $\text{sign}(\bar{f}(\mathbf{x})) = (\text{sign } f(\mathbf{x}))$. A bunch of examples are given in [6]. Note the result that the lowly squared difference $\mathcal{L}(y,f) = (y - f)^2$ leads to a large margin classifier since if $|y| = 1$, then $(y - f)^2 \equiv (1 - yf)^2$. This large margin classifier (!) is sometimes called the least squares support vector machine, but it is nothing more than ordinary ridge regression on data that have been coded as $\pm 1$. Many large margin classifiers have been proposed, both convex and nonconvex, that claim various properties; four of the many are described in [11, 14, 17, 19]. These classifiers are said to have some special advantages, either theoretical, computational or practical, and it is interesting to understand more generally the circumstances under which one cost function can be better than another. Considering accuracy as well as computational tractability, it is unlikely that there will be just one best cost function for all classification problems (see the comparison in Figure 1). The hinge function occupies a niche as a general purpose large margin classifier that is the closest convex upper bound, in some sense, to the misclassification function.

## 4. PROBABILITY ESTIMATES AND THE SVM

I respectfully disagree with the authors' remark that "from a statistical point of view, an important subject remains open: the interpretability of the SVM outputs." I think the appropriate interpretation is that the SVM targets the sign of the log odds ratio *directly*; see [7]. Since the target function sign $f(\mathbf{x})$ is discontinuous at $f(\mathbf{x})$, and the SVM is found as an optimization problem in a RKHS which is typically a space of continuous functions, it cannot jump at the boundary, but there may be a Gibbs effect there. Since the SVM is generally a smooth approximation which tends not to stray too far outside of the interval $[-1, 1]$, there is a tendency to believe that $2p - 1$ can be inferred from the SVM. This is not, however the case. A toy problem which is easy to drive toward asymptopia illustrates this point.

In Figure 2, the solid line gives $2p(x) - 1$, where $p(x)$ is the true conditional probability of the + class. Data $y_i$ have been generated as $y_i = 1$ with probability $p(x_i)$ and $-1$ with probability $1 - p(x_i)$ for 300 equally spaced points $x(i)$ in the interval $[-2, 2]$. The logit $f(x)$ has been estimated as $\hat{f}(x)$ via penalized likelihood, and $2\hat{p}(x) - 1$ is plotted as the dotted line, where $\hat{p} = e^{\hat{f}}/(1 + e^{\hat{f}})$. The dashed line gives the support vector machine estimate from the same data. It can be seen that the SVM is trying to estimate $-1$ for $x < 0$ and $+1$ for $x > 0$, which is the Bayes optimal classifier here. A small Gibbs effect near the class boundary $x = 0$ is evident, although the penalized likelihood and SVM will essentially pick out the same classification boundary. Further examples of this phenomenon in the context of the multicategory SVM of Lee, Lin and Wahba can be found in [5]. A comparative discussion of the multicategory SVM and a multicategory penalized likelihood estimate can be found in [16].

## 5. SUPPORT VECTOR REGRESSION

A precursor of the $\varepsilon$ insensitive loss function can be found in [15], where the loss function is $L(y, f(\mathbf{x}_i)) = 0$ if $|y - f| \leq \varepsilon$ and $\infty$ otherwise. In 1969 only highly quantized data were available from satellites, but computation of such estimates was iffy.

## 6. SPARSITY, VARIABLE SELECTION

In many classification problems, it is desirable to learn which components of the proposed attribute vector are actually contributing substantially to the actual classification. Two recent contributions are

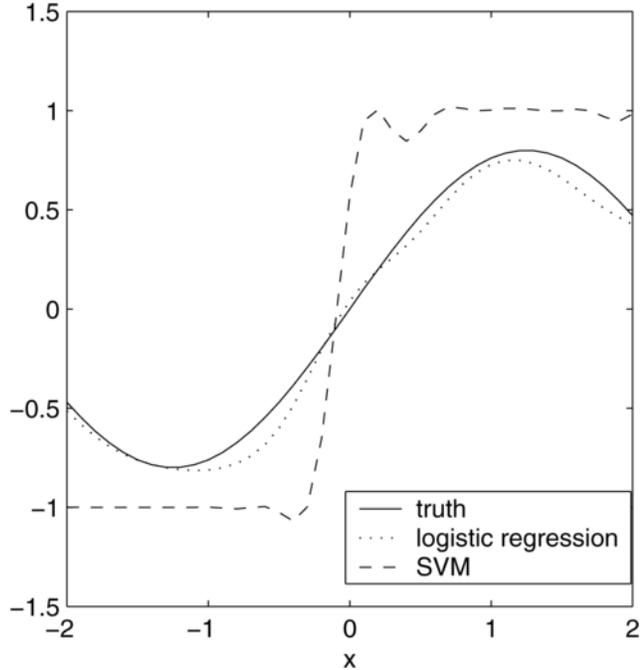

FIG. 2. *Solid line: true conditional probability $2p(x) - 1 = \text{prob}\, Y = 1$; dashed line: fitted SVM; dotted line: fitted penalized likelihood estimate. Data $y_i$ have been generated according to $p(x)$ for 300 equally spaced values of $x$.*

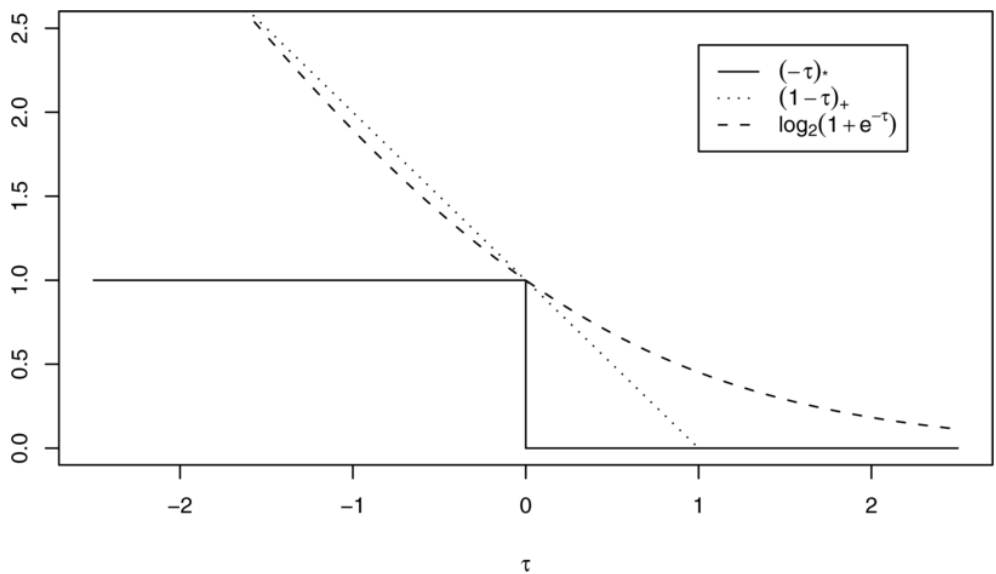

FIG. 1. *Comparison of the cost functions $c(\tau) = (-\tau)_+$, $c(\tau) = (1 - \tau)_+$ and $c(\tau) = \log_2(1 + e^{-\tau})$, which are the misclassification function, the hinge function and the negative log-likelihood function, respectively. Any strictly convex function that goes through 1 at $\tau = 0$ will be an upper bound on the misclassification function $(-\tau_+)$ and will be a looser bound than some hinge function $(1 - \theta\tau)_+$.*

[4] and [18]. The trick is to add $\ell_1$ (absolute value) penalties on coefficients of variables or terms in the penalty functional, which induces sparsity, as is well known. An early proponent of adding $\ell_1$ penalties in classification algorithms to promote sparsity is [2]; there are many other recent related contributions. In practice the major challenge in many problems involves which attributes, or clusters of attributes, to put into the model to begin with. This challenge appears in images, sounds, handwriting, text, genomic data, meteorological data, astronomical data and elsewhere. Many open questions remain in particular contexts.

## 7. REGULARIZED KERNELS FROM DISSIMILARITY DATA

Some recent work [9] focused on fitting kernels from noisy, scattered, incomplete dissimilarity data, which can then be used as a dimension reduction tool or in a SVM or multicategory SVM. Given a set of objects (protein sequences in [9]) and dissimilarity information $d_{ij}$ between the $i$th and $j$th object, for a sufficiently rich subset $\Omega$ of the $\binom{n}{2}$ pairs, one finds an $n \times n$ kernel (nonnegative definite matrix) $K_\mu$ over "object space" to yield

$$(1) \qquad \min_{K \in \mathcal{S}} \sum_{ij \in \Omega} |d_{ij} - \hat{d}_{ij}| + \mu \operatorname{trace} K,$$

where $\mathcal{S}$ is the class of nonnegative definite $n \times n$ matrices and $\hat{d}_{ij} = K_\mu(i,i) + K_\mu(j,j) - 2K_\mu(i,j)$. It is necessary to choose $\mu$ and it is useful to truncate the eigenvalues of $K_\mu$ after the first $p$, where $p$ can be chosen so as to retain some specified percentage of the trace. Suppressing $\mu$ and the truncation level $p$, a support vector machine $f(i)$, $i = 1, \ldots, n$, can be defined in object space as

$$f(i) = \sum_{\ell=1}^n c_\ell K(i,\ell)$$

by minimizing

$$\sum_{i=1}^n (1 - y_i f(i))_+ + \mu c' K^\dagger c$$

or its multicategory analog from [5].

To classify a new object ($i = n+1$), the "newbie" algorithm is used. It goes as follows: Given $d_{i,n+1}$ for sufficiently many $i$, find $b \in E^n$ and constant $c$ to minimize

$$\sum_i |d_{i,n+1} - \hat{d}_{i,n+1}|$$

over $b \in \operatorname{range}(K)$ and $c - b'K^\dagger b \geq 0$. The $b$ and $c$ are used to give a new $(n+1) \times (n+1)$ nonnegative definite matrix with $K$ in the upper left block, and $K(n+1, n+1) = c$, $K(i, n+1) = b_i, i = 1, \ldots, n$, and $\hat{d}_{ij} = K(n+1, n+1) + K(i,i) - 2K(i, n+1)$. Then the classifier evaluated at the $(n+1)$st object is

$$f(n+1) = \sum_{\ell=1}^n c_\ell K(n+1, \ell).$$

Pseudo-attribute vectors may be defined as $\mathbf{x}(i) = (\sqrt{\lambda_1}\phi_1(i), \ldots, \sqrt{\lambda_p}\phi_p(i))$, where the $\{\lambda_\nu, \phi_\nu\}$ are the eigenvalues and eigenvectors of $K$. The newbie can be placed in this pseudo-attribute coordinate system by using its fitted distance from a sufficiently large subset of the fitted training set distances. Since $K(i,j) = (\mathbf{x}(i), \mathbf{x}(j))$, the resulting SVM is linear in the pseudo-attribute vectors. However, other SVMs can be built on the labeled pseudo-attribute vectors.

The so-called semisupervised version of this problem occurs when only some of the original training objects are labeled. Thus, there are three kinds of objects: (1) those that are in the training set and labeled; (2) those that are in the training set and not labeled, but are used to determine the geometry of the object space; and (3) unlabeled newbies. Both kinds of unlabeled data can then be classified by the SVM.

The tuning parameter $\mu$ in equation (1) can be tuned by leaving out pairs of objects (CV2) and comparing their observed distances with their fitted (pseudo-attribute) distances for a range of $\mu$.

## 8. ROBUST MANIFOLD UNROLLING

A related problem occurs when the objects of interest are believed to lie in a low-dimensional (nonlinear) manifold in some higher-dimensional space. Here then it is desired to "flatten out" the manifold and reduce the dimension before carrying out a classification or regression operation. Recent references can be found in [10], where we proposed an approach related to that in equation (1) with two modifications: (a) only distances between $k$ nearest neighbors will be used and (b) $\mu \operatorname{trace} K$ is replaced by $-\mu \operatorname{trace} K$. The effect on the resulting pseudo-attribute vectors is that they tend to "flatten out" or "unroll" due to the fact that only nearest neighbor distances are used, as well as the fact that the minus sign propels distant objects to become more distant. A longer discussion of the rationale behind this algorithm and demonstrations of its behavior are found in [10]. The semisupervised version of this problem can be defined similarly, with many potential applications. Both of these optimization problems can be solved numerically using convex cone optimization code.

## 9. WHERE ARE WE GOING?

The theory, computation and application of classification problems that relate to support vector machines and other regularization based classifiers is by no means finished work, although the extent of work so far is breathtaking. Many problems remain. Using subject matter knowledge to build kernels that embody subject matter information efficiently in various fields remains an interesting challenge. For example, text and language processing have interesting problems that involve complex relationships between components of text. Huge attribute vectors and small training sets as occur in genetic data of various kinds present their own challenges, as does the merging of heterogenous kinds of information. Multiple correlated inputs and outputs provide challenges. Improved systematic ways to choose important attributes or groups of attributes remain to be found. As the authors note, the relationships between statistical learning theory based on Vapnik–Chervonenkis dimension and SVM theory based on regularization remain to be understood better, as do regularization based approaches and other approaches to classification. Collaboration between statisticians, computer scientists, mathematicians and subject matter experts will no doubt be needed for many of the practical challenges.

## ACKNOWLEDGMENTS

This research was supported by NSF Grant DMS-00-7292, by ONR Grant NN00140610 095 and by NIH Grant EY09946.